\documentclass[times,authoryear,12pt]{elsarticle}
\date{\today}
\usepackage{natbib}
\usepackage{geometry,txfonts}
\usepackage[labelfont=bf]{caption}
\usepackage{amsmath,amsfonts,amssymb,amsthm,booktabs,color,epsfig,graphicx,url}
\urldef\myurl\url{www.ualberta.ca/~dwiens/}
\usepackage{appendix}
\usepackage[utf8]{inputenc}
\usepackage[english]{babel}
 \makeatletter
\def\ps@pprintTitle{%
  \let\@oddhead\@empty
  \let\@evenhead\@empty
  \def\@oddfoot{\reset@font\hfil\thepage\hfil}
  \let\@evenfoot\@oddfoot
}
\makeatother
 
\journal{U. Alberta preprint series}

\begin{document}
\begin{frontmatter}
\title{Jittering and Clustering:\linebreak Strategies for the Construction of Robust Designs}

\author[A1]{Douglas P. Wiens\corref{mycorrespondingauthor}}

\address[A1]{Mathematical \& Statistical Sciences,
	University of Alberta,
	Edmonton, Canada,  T6G 2G1
 \newline \newline \today}

\cortext[mycorrespondingauthor]{E-mail: \url{doug.wiens@ualberta.ca}. 
Supplementary material is at {\myurl}.}

\begin{abstract}
We discuss, and give examples of, methods for
randomly implementing some minimax robust designs from the literature. These
have the advantage, over their deterministic counterparts, of having bounded
maximum loss in large and very rich neighbourhoods of the, almost certainly
inexact, response model fitted by the experimenter. Their maximum loss
rivals that of the theoretically best possible, but not implementable,
minimax designs. The procedures are then extended to more general robust
designs. For two-dimensional designs we sample from contractions of Voronoi
tessellations, generated by selected basis points, which partition the
design space. These ideas are then extended to $k$-dimensional designs for
general $k$.
\end{abstract}

\end{frontmatter}

\section{Introduction and summary \label{section: intro}}

In this article we investigate various methods of implementing experimental
designs, robust against model inadequacies. We begin with a review of the
`minimax' theory of robustness of design, and of some minimax designs from
the literature. For this we initially follow Atkinson (1996) and view a
design as any probability measure on the design space. It will be seen that
the designs which protect against a large class of alternative response
models are necessarily absolutely continuous, and so lose their optimality
when approximated by implementable, discrete (deterministic) designs. Two
remedies for this and other issues are proposed, suggested by work of Waite
and Woods (2022), who propose and study \textit{random design} strategies.

The first remedy is a random design strategy termed \textit{jittering}. The
designs are obtained by uniform sampling from small neighbourhoods of an
optimal set $t^{\ast }=\left\{ t_{i}|i=1,...,n\right\} $ of points, chosen
to represent the minimax design density. Both completely random and
stratified random -- i.e.\ random within each neighbourhood -- are
considered. We assess these designs by looking at the sample distributions
of the mean squared prediction errors incurred; with respect to these
measures both sampling strategies typically lead to designs very nearly
optimal, with the stratification strategy clearly outperforming its
completely random counterpart.

We then investigate a strategy leading to \textit{cluster} designs,
motivated by the observation that robust designs for a particular response
model tend to place their mass near those points $t_{i}\in $ $t^{\ast }$ at
which classically optimal designs, focussed solely on variance minimization,
are replicated -- but with their support points spread out in clusters of
nearby points, rather than being replicated. In clustering the idea is to
sample from densities concentrated near the $t_{i}$. An advantage to this
method over jittering is that there is no need for the minimax design to
already have been derived.

Both these approaches parallel the `random translation design strategy' of
Waite and Woods (2022), who sample uniformly in small neighbourhoods of a
chosen set of points, but with some significant differences. The choice of $%
t^{\ast }$ in jittering allows for designs whose maximum expected loss
rivals that of the minimax, absolutely continuous design. In clustering,
both the support of the non-uniform densities from which we sample, and the
extent of their concentration near the $t_{i}$, are governed by a
user-chosen parameter $\nu $, representing the bias/variance trade-off
desired by a user seeking robustness against model misspecifications.

We start by applying these ideas in several one-dimensional cases for which minimax designs -- in continuous or discrete design spaces -- have previously been derived. The
framework is that the experimenter will fit a polynomial response, and our
random designs have points assigned at random but in a structured manner
near the $t^{\ast }$. The densities from which we sample are chosen to
capture the salient properties of the minimax designs (in jittering) or classically
optimal but deterministic designs (in clustering). The structure we impose
-- especially that of stratification -- is shown, through a number of
examples, to lead to efficient designs approximating the variance minimizing
properties of the deterministic designs concentrated on the $t^{\ast }$. But
the randomness, leading to the clustering effect -- this alone is known to
increase robustness -- ensures that the bias is bounded as well, even in
continuous design spaces in which the bias of deterministic designs can be
unbounded. 

We then consider two-dimensional clustering applications in which intervals
containing the $t_{i}$\ are replaced by less regular regions formed by
shrinking Voronoi tessellations generated by $t^{\ast }$. We sample from
spherical beta densities centred on the $t_{i}$, and suggest tuning
constants which again result in both efficiency and robustness. We finish
with recommendations for the construction of $k$-dimensional designs for $%
k\geq 3$. 

The examples were prepared using \textsc{matlab}; the code is available on
the author's website.

\section{Minimax robustness of design \label{section: minimax}}

The theory of robustness of design was largely initiated by Box and Draper
(1959), who investigated the robustness of some classical experimental
designs in the presence of certain model inadequacies, e.g.\ designs optimal
for a low order polynomial response when the true response was a polynomial
of higher order. Huber (1975) derived \textit{minimax} designs for straight
line regression; these minimize the maximum integrated mean squared error,
with the maximum taken over a large class of alternative responses. Wiens
(1990, 1992) extended these results to multiple regression responses and in
a variety of other directions -- see Wiens (2015) for a summary of these and
other approaches to robustness of design. Specifically, the general problem
is phrased in terms of an approximate regression response 
\begin{equation}
E\left[ Y\left( \boldsymbol{x}\right) \right] \approx \boldsymbol{f}^{\prime
}\left( \boldsymbol{x}\right) \boldsymbol{\theta },  \label{approx}
\end{equation}%
for $p$ regressors $\boldsymbol{f}$, each functions of $q$ independent
variables $\boldsymbol{x}$, and a parameter $\boldsymbol{\theta }$. Since (%
\ref{approx}) is an approximation the interpretation of $\boldsymbol{\theta }
$ is unclear; we \textit{define }this target parameter by 
\begin{equation}
\boldsymbol{\theta }=\arg \min_{\boldsymbol{\eta }}\int_{\mathcal{X}}\left( E%
\left[ Y\left( \boldsymbol{x}\right) \right] -\boldsymbol{f}^{\prime }\left( 
\boldsymbol{x}\right) \boldsymbol{\eta }\right) ^{2}\mu \left( d\boldsymbol{x%
}\right) ,  \label{theta def}
\end{equation}%
where $\mu \left( d\boldsymbol{x}\right) $ represents either Lebesgue
measure or counting measure, depending upon the nature of the \textit{design
space} $\mathcal{X}$. We then define $\psi \left( \boldsymbol{x}\right) =E%
\left[ Y\left( \boldsymbol{x}\right) \right] -\boldsymbol{f}^{\prime }\left( 
\boldsymbol{x}\right) \boldsymbol{\theta }$. This results in the class of
responses $E\left[ Y\left( \boldsymbol{x}\right) \right] =\boldsymbol{f}%
^{\prime }\left( \boldsymbol{x}\right) \boldsymbol{\theta }+\psi \left( 
\boldsymbol{x}\right) $, with -- by virtue of (\ref{theta def}) -- $\psi $
satisfying the orthogonality requirement 
\begin{equation}
\int_{\mathcal{X}}\boldsymbol{f}\left( \boldsymbol{x}\right) \psi \left( 
\boldsymbol{x}\right) \mu \left( d\boldsymbol{x}\right) =\boldsymbol{0}.\ 
\label{orthogonality}
\end{equation}%
Assuming that $\mathcal{X}$ is rich enough that the matrix\ $\boldsymbol{A}%
=\int_{\mathcal{X}}\boldsymbol{f}\left( \boldsymbol{x}\right) \boldsymbol{f}%
^{\prime }\left( \boldsymbol{x}\right) \mu \left( d\boldsymbol{x}\right) $
is invertible, the parameter defined by (\ref{theta def}) and (\ref%
{orthogonality}) is unique.

We identify a design with its design measure -- a probability measure $\xi
\left( d\boldsymbol{x}\right) $ on $\mathcal{X}$. Define 
\begin{equation*}
\boldsymbol{M}_{\xi }=\int_{\mathcal{X}}\boldsymbol{f}\left( \boldsymbol{x}%
\right) \boldsymbol{f}^{\prime }\left( \boldsymbol{x}\right) \xi \left( d%
\boldsymbol{x}\right) ,\text{ \ }\boldsymbol{b}_{\psi ,\xi }=\int_{\mathcal{X%
}}\boldsymbol{f}\left( \boldsymbol{x}\right) \psi \left( \boldsymbol{x}%
\right) \xi \left( d\boldsymbol{x}\right) ,
\end{equation*}%
and assume $\xi $ is such that $\boldsymbol{M}_{\xi }$ is invertible. The
covariance matrix of the least squares estimator $\boldsymbol{\hat{\theta}}$%
, assuming homoscedastic errors with variance $\sigma _{\varepsilon }^{2}$,
is $\left( \sigma _{\varepsilon }^{2}/n\right) \boldsymbol{M}_{\xi }^{-1}$,
and the bias is $E\left[ \boldsymbol{\hat{\theta}-\theta }\right] =%
\boldsymbol{M}_{\xi }^{-1}\boldsymbol{b}_{\psi ,\xi }$; together these yield
the mean squared error (\textit{mse}) matrix 
\begin{equation*}
\text{\textsc{mse}}\left[ \boldsymbol{\hat{\theta}}\right] =\frac{\sigma
_{\varepsilon }^{2}}{n}\boldsymbol{M}_{\xi }^{-1}+\boldsymbol{M}_{\xi }^{-1}%
\boldsymbol{b}_{\psi ,\xi }\boldsymbol{b}_{\psi ,\xi }^{\prime }\boldsymbol{M%
}_{\xi }^{-1}
\end{equation*}%
of the parameter estimates, whence the \textit{mse} of the fitted values $%
\hat{Y}\left( \boldsymbol{x}\right) =\boldsymbol{f}^{\prime }\left( 
\boldsymbol{x}\right) \boldsymbol{\hat{\theta}}$ is%
\begin{equation*}
\text{\textsc{mse}}\left[ \hat{Y}\left( \boldsymbol{x}\right) \right] =\frac{%
\sigma _{\varepsilon }^{2}}{n}\boldsymbol{f}^{\prime }\left( \boldsymbol{x}%
\right) \boldsymbol{M}_{\xi }^{-1}\boldsymbol{f}\left( \boldsymbol{x}\right)
+\left( \boldsymbol{f}^{\prime }\left( \boldsymbol{x}\right) \boldsymbol{M}%
_{\xi }^{-1}\boldsymbol{b}_{\psi ,\xi }\right) ^{2}.
\end{equation*}%
\ 

\noindent A loss function that is commonly employed is the \textit{%
integrated mse} of the predictions: 
\begin{eqnarray}
\text{\textsc{imse}}\left( \xi |\psi \right) &=&\int_{\mathcal{X}}\text{%
\textsc{mse}}\left[ \hat{Y}\left( \boldsymbol{x}\right) \right] d\boldsymbol{%
x}  \notag \\
&=&\frac{\sigma _{\varepsilon }^{2}}{n}tr\left( \boldsymbol{AM}_{\xi
}^{-1}\right) +\boldsymbol{b}_{\psi ,\xi }^{\boldsymbol{\prime }}\boldsymbol{%
M}_{\xi }^{-1}\boldsymbol{AM}_{\xi }^{-1}\boldsymbol{b}_{\psi ,\xi }+\int_{%
\mathcal{X}}\psi ^{2}\left( \boldsymbol{x}\right) \mu \left( d\boldsymbol{x}%
\right) \mathbf{.}  \label{imse}
\end{eqnarray}%
The dependence on $\psi $ is eliminated by adopting a \textit{minimax}
approach, according to which one first maximizes (\ref{imse}) over a
neighbourhood of the assumed response. This neighbourhood is constrained by (%
\ref{orthogonality}) and by a bound $\int_{\mathcal{X}}\psi ^{2}\left( 
\boldsymbol{x}\right) \mu \left( d\boldsymbol{x}\right) \leq \tau ^{2}/n$,
required so that errors due to bias and to variation remain of the same
order, asymptotically.

Huber (1975) took $\mathcal{X}$ to be an interval of the real line and
assumed that the minimax design measure had a density $m\left( x\right) $;
Wiens (1992) justified this assumption by proving that any design whose
design space $\mathcal{X}$ has positive Lebesgue measure, and which places
positive mass on a set of Lebesgue measure zero, necessarily has $\sup_{\psi
}$\textsc{imse}$\left( \xi |\psi \right) =\infty $. Thus in order that a
design on an interval, hypercube, etc. have finite maximum loss, it must be
absolutely continuous. For such a design $\max_{\psi }$\textsc{imse}$\left(
\xi |\psi \right) $ is $\left( \sigma _{\varepsilon }^{2}+\tau ^{2}\right)
/n $ times 
\begin{equation}
I_{\nu }\left( \xi \right) =\left( 1-\nu \right) tr\boldsymbol{AM}_{\xi
}^{-1}+\nu ch_{\max }\boldsymbol{K}_{\xi }\boldsymbol{H}_{\xi }^{-1},
\label{max loss}
\end{equation}%
where 
\begin{equation*}
\boldsymbol{H}\mathbf{_{\xi }}=\boldsymbol{M}\mathbf{_{\xi }}\boldsymbol{A}%
^{-1}\boldsymbol{M}_{\xi },\text{ }\boldsymbol{K}_{\xi }=\int_{\mathcal{X}}%
\boldsymbol{f}\left( \boldsymbol{x}\right) \boldsymbol{f}^{\prime }\left( 
\boldsymbol{x}\right) m^{2}\left( \boldsymbol{x}\right) d\boldsymbol{x},
\end{equation*}%
$ch_{\max }$ denotes the maximum eigenvalue and $\nu =\tau ^{2}\left/ \left(
\sigma _{\varepsilon }^{2}+\tau ^{2}\right) \right. \in \left[ 0,1\right] $,
representing the relative importance, to the experimenter, of errors due to
bias rather than to variance. Our examples in this article use $\nu = .5$; other values tell much the same story.

With 
\begin{eqnarray}
\boldsymbol{G}_{\xi } &=&\boldsymbol{K}_{\xi }-\boldsymbol{H}_{\xi }  \notag
\\
&=&\int_{\mathcal{X}}\left[ \left( m\left( \boldsymbol{x}\right) \boldsymbol{%
I}_{p}-\boldsymbol{M}_{\xi }\boldsymbol{A}^{-1}\right) \boldsymbol{f}\left( 
\boldsymbol{x}\right) \right] \left[ \left( m\left( \boldsymbol{x}\right) 
\boldsymbol{I}_{p}-\boldsymbol{M}_{\xi }\boldsymbol{A}^{-1}\right) 
\boldsymbol{f}\left( \boldsymbol{x}\right) \right] ^{\boldsymbol{\prime }}d%
\boldsymbol{x}\mathbf{,}  \notag \\
\mathbf{r}_{\xi }\left( \boldsymbol{x}\right) &=&\frac{\tau }{\sqrt{n}}%
\boldsymbol{G}_{\xi }^{-1/2}\left( m\left( \boldsymbol{x}\right) \boldsymbol{%
I}_{p}-\boldsymbol{M}_{\xi }\boldsymbol{A}^{-1}\right) \boldsymbol{f}\left( 
\boldsymbol{x}\right) ,  \label{r}
\end{eqnarray}%
the least favourable contaminant is%
\begin{equation}
\psi _{\xi }\left( \boldsymbol{x}\right) =\mathbf{r}_{\xi }^{\prime }\left( 
\boldsymbol{x}\right) \mathbf{\beta _{\xi },\ }  \label{least favourable}
\end{equation}%
where $\mathbf{\beta }_{\xi }$ is the unit eigenvector belonging to the
maximum eigenvalue of $\boldsymbol{G}_{\xi }^{1/2}\boldsymbol{H}_{\xi }^{-1}%
\boldsymbol{G}_{\xi }^{1/2}+\boldsymbol{I}_{p}$. See Wiens (2015) for
details and further references.

\subsection{Random designs\label{section: random designs}}

In the following sections we construct distributions $\Phi \left( 
\boldsymbol{x}\right) $, with densities $\phi \left( \boldsymbol{x}\right) $%
, and propose randomly choosing design points from $\Phi $. An $n$-point
design $D=\left\{ \boldsymbol{x}_{i}\right\} _{i=1}^{n}$ chosen in this way
has design measure $\delta =n^{-1}\sum \delta _{\boldsymbol{x}_{i}}$, where $%
\delta _{\boldsymbol{x}_{i}}$ is point mass at $\boldsymbol{x}_{i}\sim \Phi $%
. By the preceding any such design has unbounded \textsc{imse} once it is
chosen. Of interest however is the \textit{expected} \textsc{imse} against a
common alternative $\psi $; for this we take the least favourable
contaminant $\psi _{\Phi }$, given by (\ref{r}) and (\ref{least favourable})
but with $\xi $ replaced by $\Phi $. In the Appendix we show that 
\begin{equation}
E_{\Phi }\left[ \text{\textsc{imse}}\left( \delta |\psi _{\Phi }\right) %
\right] =\left( \sigma _{\varepsilon }^{2}+\tau ^{2}\right) /n\times J_{\nu
}\left( \Phi \right) ,  \label{E(imse)}
\end{equation}%
where%
\begin{eqnarray}
J_{\nu }\left( \Phi \right) &=&E_{\Phi }\left[ j_{\nu }\left( \delta \right) %
\right] \text{, for }  \label{J} \\
j_{\nu }\left( \delta \right) &=&\left( 1-\nu \right) tr\boldsymbol{AM}%
_{\delta }^{-1}+\nu \gamma _{\delta }\text{, and}  \notag \\
\gamma _{\delta } &=&\mathbf{\beta }_{\Phi }^{\prime }\boldsymbol{G}_{\Phi
}^{-1/2}\left( \boldsymbol{M}_{\phi }\boldsymbol{M}_{\delta }^{-1}-%
\boldsymbol{M}_{\Phi }\boldsymbol{A}^{-1}\right) \boldsymbol{A}\left( 
\boldsymbol{M}_{\delta }^{-1}\boldsymbol{M}_{\phi }-\boldsymbol{A}^{-1}~%
\boldsymbol{M}_{\Phi }\right) \boldsymbol{G}_{\Phi }^{-1/2}\mathbf{\beta }%
_{\Phi }+1\text{.}  \notag
\end{eqnarray}%
Here $\boldsymbol{M}_{\delta }\overset{def}{=}\frac{1}{n}\sum_{\boldsymbol{x}%
_{i}\in D}\boldsymbol{f}\left( \boldsymbol{x}_{i}\right) \boldsymbol{f}%
^{\prime }\left( \boldsymbol{x}_{i}\right) $ and $\boldsymbol{M}_{\phi }%
\overset{def}{=}\frac{1}{n}\sum_{\boldsymbol{x}_{i}\in D}\boldsymbol{f}%
\left( \boldsymbol{x}_{i}\right) \phi \left( \boldsymbol{x}_{i}\right) 
\boldsymbol{f}^{\prime }\left( \boldsymbol{x}_{i}\right) $; $\mathbf{\beta }%
_{\Phi }$ is the unit eigenvector belonging to the maximum eigenvalue of $%
\boldsymbol{G}_{\Phi }^{1/2}\boldsymbol{H}_{\Phi }^{-1}\boldsymbol{G}_{\Phi
}^{1/2}+\boldsymbol{I}_{p}$. \ 

Note that both $\boldsymbol{M}_{\delta }$ and $\boldsymbol{M}_{\phi }$ are
random. The expectation in (\ref{J}) can be estimated by averaging over a
large number of realizations of $\delta $ -- we do this in Sections \ref%
{section: clustering1d} and \ref{section: clusteringkd}. In the special case
that $\phi \left( \boldsymbol{x}\right) $ is constant on its support -- as
is the case in \S \ref{section: jittering} -- $\boldsymbol{M}_{\delta }^{-1}%
\boldsymbol{M}_{\phi }$ is a constant multiple of $\boldsymbol{I}_{p}$, $%
\gamma _{\delta }$ is non-random, and these formulas simplify considerably
-- see (\ref{easy J}).

An efficient design strategy should result in $J_{\nu }\left( \Phi \right) $
being close to $I_{\nu }\left( \Phi \right) $, with the $j_{\nu }\left(
\delta \right) $ being concentrated near their expectation.

A referee has pointed out that a more natural measure is perhaps the
maximizer $\psi _{0}$ of $E_{\delta }\left[ \text{\textsc{imse}}\left(
\delta |\psi \right) \right] $; this turns out to be computationally
infeasible in all but the simplest scenarios. And see \S \ref{section: SLR},
where we argue that a contaminant less favourable than $\psi _{\Phi }$ is
difficult to imagine. See also Figure \ref{fig:cluster1}(d)-(f).

\section{Jittering \label{section: jittering}}

There are obvious issues in implementing an absolutely continuous design
measure within this framework, since any discrete approximation necessarily
suffers from the drawback, as above, that the maximum loss is infinite.
Noting that in this case the least favourable contaminating function $\psi $
is largely concentrated on a set of measure zero -- an unlikely eventuality
against which to seek protection -- Wiens (1992, p.\ 355) states that
\textquotedblleft Our attitude is that an approximation to a design which is
robust against more realistic alternatives is preferable to an exact
solution in a neighbourhood which is unrealistically
sparse.\textquotedblright\ He places one observation at each of the
quantiles 
\begin{equation}
t_{i}=\xi ^{-1}\left( \frac{i-1/2}{n}\right) ,\ i=1,...,n,  \label{quant}
\end{equation}%
which is the $n$-point design closest to $\xi $ in Kolmogorov distance (Fang
and Wang 1994; see Xu and Yuen 2011 for other possibilities).

Despite the disclaimer above, such discrete implementations have become
controversial; see in particular Bischoff (2010). In this article we
investigate a resolution to these difficulties offered by Waite and Woods
(2022), who propose randomly sampling the design points from uniform
densities highly concentrated in small neighbourhoods of an optimally chosen
set of deterministic points. In our case we propose random sampling from a
piecewise uniform density 
\begin{equation}
\phi _{n}\left( x;c\right) =\frac{1}{2c}\sum_{i=1}^{n}I\left[ t_{i}-\frac{c}{%
n}\leq x\leq t_{i}+\frac{c}{n}\right] ,  \label{phi}
\end{equation}%
for chosen $c\in \left( 0,1\right) $.

We illustrate the method in the context of straight line regression -- $%
\mathcal{X}=\left[ -1,1\right] $ and $\boldsymbol{f}\left( x\right) =\left(
1,x\right) ^{\prime }$ -- for which Huber (1975) obtained the minimax density%
\begin{equation*}
m\left( x\right) =3\left( x^{2}-\alpha \right) ^{+}/d\left( \alpha \right) ,
\end{equation*}%
with $\alpha $ chosen to minimize (\ref{max loss}), which in terms of 
\begin{equation*}
\mu _{2}(\alpha )=\int_{-1}^{1}x^{2}m\left( x\right) dx,\text{ }\kappa
_{0}(\alpha )=\int_{-1}^{1}m^{2}\left( x\right) dx,\text{ }\kappa
_{2}(\alpha )=\int_{-1}^{1}x^{2}m^{2}\left( x\right) dx,
\end{equation*}%
is%
\begin{equation*}
K_{\nu }(\alpha )=2\left( 1-\nu \right) \left( 1+\frac{1}{3\mu _{2}}\right)
+2\nu \max \left( \kappa _{0},\frac{\kappa _{2}}{3\mu _{2}^{2}}\right) .
\end{equation*}%
Apart from minor modifications resulting from the change in the support to $%
\left[ -1,1\right] $ from $\left[ -1/2,1/2\right] $, the details of the
construction of $m$ are as in Huber (1975). We assume that $\max \left(
\kappa _{0},\kappa _{2}/3\mu _{2}^{2}\right) =$ $\kappa _{0}$ and check this
once $m$ is obtained. We find%
\begin{equation*}
d\left( \alpha \right) =\left\{ 
\begin{tabular}{cc}
$2\left( 1-3\alpha \right) ,$ & $\alpha \leq 0,$ \\ 
$2\left( 1-\sqrt{\alpha }\right) ^{2}\left( 1+2\sqrt{\alpha }\right) ,$ & $%
\alpha \geq 0,$%
\end{tabular}%
\right.
\end{equation*}%
with $\alpha $ and $\nu $ related by%
\begin{equation*}
\nu ^{-1}=\left\{ 
\begin{tabular}{cc}
$1+\frac{9\left( 3-5\alpha \right) ^{2}}{25\left( 1-3\alpha \right) ^{3}},$
& $\alpha \leq 0,$ \\ 
$1+\frac{9\left( 3+6\sqrt{\alpha }+4\alpha +2\alpha ^{3/2}\right) ^{2}}{%
25\left( 1-\sqrt{\alpha }\right) ^{2}\left( 1+2\sqrt{\alpha }\right) ^{3}},$
& $\alpha \geq 0.$%
\end{tabular}%
\right.
\end{equation*}%
The limiting cases are (i) $\alpha \rightarrow -\infty $, $\nu \rightarrow 1$%
, $m\left( x\right) \rightarrow .5$ (the uniform density), (ii) $\alpha =0$, 
$\nu =25/106$, $m\left( x\right) =3x^{2}/2$, and (iii) $\alpha \rightarrow
\infty $, $\nu \rightarrow 0$, $m\left( x\right) \rightarrow $ point masses
of $1/2$ at $\pm 1$.

It is a fortuitous consequence of the choice of \textit{imse }as loss that
for all $\nu \in \left[ 0,1\right] $, $\max \left( \kappa _{0},\kappa
_{2}/3\mu _{2}^{2}\right) =$ $\kappa _{0}$, the choice used in the
derivation of the minimizing density $m$. For other common choices -- D-, A-
and E-optimality for instance -- the situation is far more complicated. See
Daemi and Wiens (2013).

\subsection{Jittered designs for SLR\label{section: SLR}}

In the construction of the sampling density (\ref{phi}) for this example we
will take $\alpha \leq 0$ -- the case of most interest from a robustness
standpoint -- and then for $m$ as above, the symmetrically placed points $%
t_{i}$ are determined by%
\begin{equation*}
t_{i}^{3}-3\alpha t_{i}=\left( 1-3\alpha \right) \left( \frac{2i-1-n}{n}%
\right) ,\ i=1,...,n.
\end{equation*}%
This equation has an explicit solution furnished by Cardano's formula:%
\begin{equation*}
t_{i}=\left( -s/2+\sqrt{\Delta }\right) ^{1/3}+\left( -s/2-\sqrt{\Delta }%
\right) ^{1/3},
\end{equation*}%
for 
\begin{equation*}
s=-\left( 1-3\alpha \right) \left( \frac{2i-1-n}{n}\right) ,\text{ }\Delta
=s^{2}/4-\alpha ^{3}>0.
\end{equation*}

From (\ref{quant}), and the bowl-shape of $m(x)$, one infers that the
distances between adjacent $t_{i}$ are smallest near $\pm 1$, largest near $0
$. Thus the intervals of support of $\phi_n $ will be non-overlapping, and
within $\left[ -1,1\right] $, as long as $t_{1}-c/n\geq -1$, i.e.\ $c\leq
n\left( 1+t_{1}\right) $. Note that the interpretation of $c$ is that it is
the proportion of the design space being randomly sampled. 

\begin{figure}[tb]
\centering
\includegraphics[scale=1]{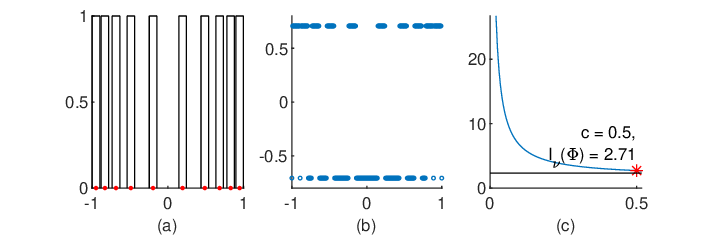}
\caption{(a) Jittered design density $\protect\phi %
_{n}\left( x;c\right) $ for an approximately linear univariate response
using $\protect\nu =.5$, $c=.5$, $n=10$. (b) $\protect\sqrt{n}\protect\psi %
_{\Phi }\left( x;c\right) /\protect\tau $. (c) $I_{\protect\nu }\left( \Phi
\right) $ vs. $c$; horizontal line at $I_{\protect\nu }\left( \protect\xi %
\right) =2.31$.}
\label{fig:jitter1}
\end{figure}

\begin{figure}[tb]
\centering
\includegraphics[scale=1]{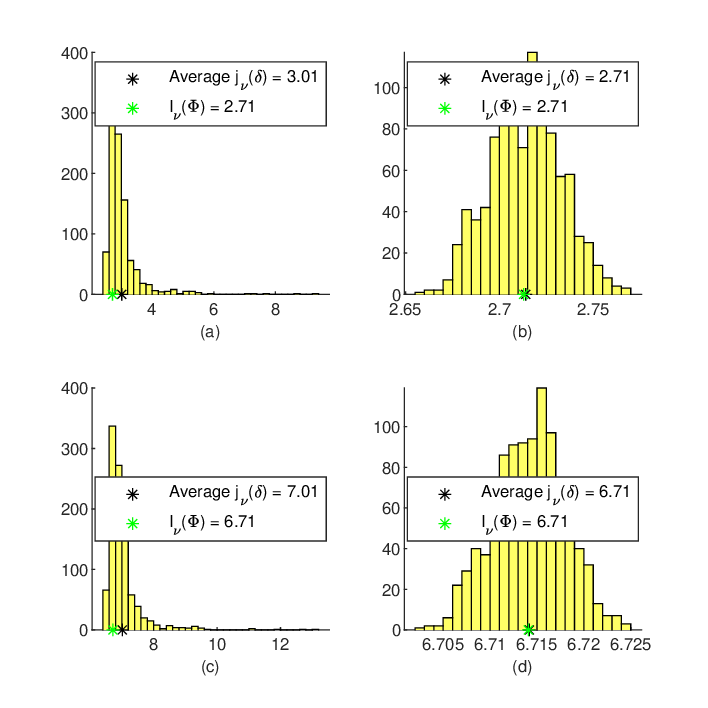}
\caption{Values of $j_{\protect\nu %
}\left( \protect\delta \right) $ and their averages, estimating $J_{\protect%
\nu }\left( \Phi \right) $, from 1000 jittered designs $\protect%
\delta $ using $\protect\nu =.5$. Top: $c=.5$; bottom: $c=.1$. 
(a), (c): Completely random sampling. (b), (d): Stratified random
sampling.}
\label{fig:jitter2}
\end{figure}

In Figure \ref{fig:jitter1}(a) we plot $\phi _{n}\left( x;c\right) $, when
placing equal weight on protection against bias versus variance ($\nu =.5$), 
$50\%$ of the design space to be sampled from $\Phi $ ($c=.5$) and $n=10$.
The $t_{i}$ are the quantiles arising from $m(x)\propto \left(
x^{2}+.325\right) $.

A comparison of the maximum loss (\ref{max loss}) of $\xi $ versus that of
the design measure $\Phi $ corresponding to $\phi $ is obtained from%
\begin{eqnarray*}
I_{\nu }\left( \xi \right) &=&2\left( 1-\nu \right) \left( 1+\frac{1}{3\mu
_{2}(\alpha )}\right) +\nu \left( 1+\frac{5}{4}\left( 3\mu _{2}(\alpha
)-1\right) ^{2}\right) , \\
I_{\nu }\left( \Phi \right) &=&2\left( 1-\nu \right) \left( 1+\frac{1}{%
3\lambda _{2}(c)}\right) +\frac{\nu }{c}\max \left( 1,\frac{1}{3\lambda
_{2}(c)}\right) ,
\end{eqnarray*}%
where 
\begin{equation*}
\mu _{2}(\alpha )=\frac{3-5\alpha }{5\left( 1-3\alpha \right) },\text{ and }%
\lambda _{2}(c)=\int_{-1}^{1}x^{2}\phi _{n}\left( x;c\right) dx=\frac{1}{n}%
\sum_{i=1}^{n}t_{i}^{2}+\frac{c^{2}}{3n^{2}}.
\end{equation*}%
As noted in \S \ref{section: random designs}, $E_{\Phi }\left[ \text{\textsc{%
imse}}\left( \delta |\psi _{\Phi }\right) \right] $ simplifies considerably
for these jittered designs and then $J_{\nu }\left( \Phi \right) $ is very
similar to $I_{\nu }\left( \Phi \right) $, plotted in Figure \ref%
{fig:jitter1}(c). We show in the Appendix that in this case (\ref{J}) becomes 
\begin{equation}
J_{\nu }\left( \Phi \right) =\left( 1-\nu \right) E_{\Phi }\left[ tr%
\boldsymbol{AM}_{\delta }^{-1}\right] +\nu \gamma _{0},\text{where }\gamma
_{0}=\frac{1}{c}\max \left( 1,\frac{1}{3\lambda _{2}\left( c\right) }\right)
,  \label{easy J}
\end{equation}%
and that, with\ $I_{S}\left( x;c\right) =I\left( \phi _{n}\left( x;c\right)
>0\right) $, the least favourable contaminant for $\Phi $ is 
\begin{equation}
\psi _{\Phi }\left( x;c\right) =\frac{\tau }{\sqrt{n}}\left( \frac{%
I_{S}\left( x;c\right) -\frac{1}{\gamma _{0}}}{\sqrt{2c\left( 1-\frac{1}{%
\gamma _{0}}\right) }}\right) \cdot \left( \frac{x}{\sqrt{\lambda _{2}}}%
\right) ^{I\left( \lambda _{2}\left( c\right) <1/3\right) }.  \label{PsiPhi}
\end{equation}%
In Figure \ref{fig:jitter1}(b) we plot a scaled version of $\psi _{\Phi
}\left( x;c\right) $. The contaminant $\psi _{\Phi }$ has the effect of
changing the uncontaminated response $E\left[ Y\left( x\right) \right]
=\theta _{0}+\theta _{1}x$ to $\left( \theta _{0}-k\right) +\theta
_{1}x+2kI\left( \phi _{n}\left( x;c\right) >0\right) $ for (when $c=.5$) $%
k=\tau /\sqrt{2n}$. Thus it biases the intercept and then places
contamination uniformly on the support of $\phi _{n}$. In the parlance of
game theory, it is difficult to see how Nature, knowing $\Phi $ but not $%
\delta $ and assumed malevolent, could respond less favourably than this.

For ease in the estimation of $J_{\nu }\left( \Phi \right) $ we note that $%
E_{\Phi }\left[ tr\boldsymbol{AM}_{\delta }^{-1}\right] =2\left\{ 1+E_{\Phi }%
\left[ \left( \mu _{\delta }^{2}+1/3\right) /\sigma _{\delta }^{2}\right]
\right\} $, where $\mu _{\delta }$ and $\sigma _{\delta }^{2}$ are the mean
and variance of the design. See Figure \ref{fig:jitter2} for comparative
values illustrating the close agreement between $J_{\nu }\left( \Phi \right) 
$ and $I_{\nu }\left( \Phi \right) $.

The plots reveal that the loss associated with the design $\Phi $ decreases
with $c$, i.e.\ as the design becomes closer to the uniform design on all of 
$\chi $, for which the bias vanishes. This is in line with the remark of Box
and Draper (1959): \textquotedblleft The optimal design in typical
situations in which both variance and bias occur is very nearly the same as
would be obtained if \textit{variance were ignored completely} and the
experiment designed so as to \textit{minimize bias alone}.\textquotedblright

\subsubsection{Sampling methods}

We constructed $1000$ completely random and stratified random designs, in
order to assess their performance. A completely random design $\delta $
consisted of $n=10$ points chosen from $\phi _{n}\left( x;c\right) $. The
resulting values of $j_{\nu }\left( \delta \right) $ are plotted in Figure %
\ref{fig:jitter2}(a),(c). Stratification consisted of choosing one design
point at random from each bin -- Figure \ref{fig:jitter2}(b),(d). The sample
averages of the losses from the randomized designs were smaller and closer
to $I_{\nu }\left( \Phi \right) $ under the stratified sampling scheme, and
more concentrated around their expectation of $J_{\nu }\left( \Phi \right) $%
, as exhibited by the much shorter tail in (b). In a further simulation, for
which the output is not displayed here, we estimated $J_{\nu }\left( \xi
\right) $, as at (\ref{J}), by drawing $1000$ samples from the minimax
density $m\left( x\right) $ and averaging their \textsc{imse}. The values $%
\left\{ j_{\nu }\right\} $ showed more variation than those plotted in
Figure \ref{fig:jitter2}(a), and with an average of $2.72$ -- significantly
larger than the target value $I_{\nu }\left( \xi \right) =2.31$. From this
we infer that jittering combined with stratification gives an efficient,
structured implementation of the minimax solution.

Simulations using other inputs also resulted in these same conclusions\ --
that our random design strategies typically yield designs very close to
optimal with respect to our robustness and efficiency requirements, and that
do not suffer from the drawback of deterministic designs of having infinite
maximum loss.

\section{Cluster designs in one dimension \label{section: clustering1d}}

Working in discrete design spaces, Wiens (2018) obtained minimax robust
designs for a variety of approximate responses. Those shown in Figure \ref%
{fig:cubic} are for cubic regression. \ The classically I-optimal design ($%
\nu =0$) minimizing integrated variance alone was derived by Studden (1977)
and places masses of $.1545$ and $.3455$ at $\pm 1$ and $\pm .4472$. The
robust designs can thus be described as taking the replicates of the
classical design and spreading their mass out (`clustering') over nearby
regions. This same phenomenon has frequently been noticed in other
situations (Fang and Wiens 2000, Heo et al. 2001 for instance).

\begin{figure}[tb]
\centering
\includegraphics[scale=1]{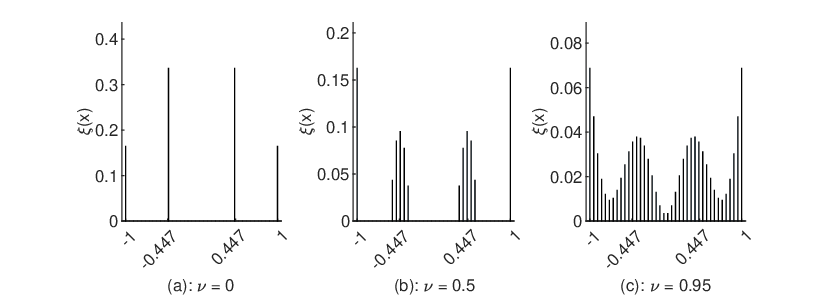}
\caption{Minimax designs for approximate
cubic regression}
\label{fig:cubic}
\end{figure}

In this section we aim to formalize this notion in order to obtain designs
competing with the minimax designs, but with finite maximum loss even in
continuous design spaces, and having the advantage of being much more easily
derived -- there is no need for the minimax designs to be known. We consider
only one-dimensional designs in this section, and will illustrate the
methods in polynomial response models of degrees $p-1=1,2,3$.

Suppose that a given static design has $p$ support points $t_{1}<\cdot \cdot
\cdot <t_{p}$ in $[-1,1]$. Define midpoints $s_{i}=\left(
t_{i}+t_{i+1}\right) /2$, $i=1,...,p-1$. Put $s_{0}=\min $ $\left(
-1,t_{1}\right) $ and $s_{p}=\max \left( 1,t_{p}\right) $. Then the $p$
intervals $I_{i}=\left\{ \left[ s_{i-1},s_{i}\right] ,\text{ }%
i=1,...,p\right\} $ cover $[-1,1]$ and have the properties that $t_{i}\in
I_{i}$ and that any point in $I_{i}$ is closer to $t_{i}$ than to any $t_{j}$%
, $j\neq i$. This is then a trivial example of a \textit{Voronoi tessellation%
}, to be considered when we pass to higher dimensions.

We propose designs consisting of points sampled from Beta densities on
subintervals of the $I_{i}$. Specifically, for $i=1,...,p$ let $c=c\left(
\nu \right) \in \left[ 0,1\right] $ satisfy $c\left( 0\right) =0$ and $%
c\left( 1\right) =1$. Put $J_{i}\left( c\right) =\left[ t_{i}-c\left(
t_{i}-s_{i-1}\right) ,t_{i}+c\left( s_{i}-t_{i}\right) \right] \equiv \left[
k_{i},l_{i}\right] $, with length $\left\vert J_{i}\right\vert =\left(
l_{i}-k_{i}\right) =c\left( s_{i}-s_{i-1}\right) =c$ $\times \left\vert
I_{i}\right\vert $. Let $\beta _{a,b}\left( x\right) $ be the Beta$\left(
a,b\right) $ density on $\left[ 0,1\right] $. Then 
\begin{equation}
\frac{1}{\left\vert J_{i}\right\vert }\beta _{a,b}\left( \frac{x-k_{i}}{%
\left\vert J_{i}\right\vert }\right) ,\text{ }x\in J_{i}\left( c\right)
\label{density}
\end{equation}%
is this density, translated and scaled to$\ J_{i}\left( c\right) $. The
interpretation of `$c$' is as before -- it is the fraction of the design
space to be sampled. Here and in the following examples we use $c=\nu ^{k}$
where $k$ is the dimension of $\boldsymbol{x}$, so that $c$ varies at the
same rate as the volume of $\chi $ as the dimensionality changes.

The parameters $\left( a_{i},b_{i}\right) $ are chosen so that the mode of (%
\ref{density}) is at $t_{i}\in J_{i}\left( c\right) $, hence the mode $%
\delta _{i}\in \left[ 0,1\right] $ of $\beta _{a,b}\left( x\right) $ is
given by 
\begin{equation*}
\delta _{i}\equiv \frac{t_{i}-k_{i}}{l_{i}-k_{i}}=\left\{ 
\begin{tabular}{cc}
$\frac{a_{i}-1}{a_{i}-1+b_{i}-1},$ & $a_{i},b_{i}>1,$ \\ 
$0,$ & $a_{i}\leq 1<b_{i},$ \\ 
$1,$ & $b_{i}\leq 1<a_{i}.$%
\end{tabular}%
\right.
\end{equation*}%
Then 
\begin{equation}
\left( a_{i}-1\right) \left( 1-\delta _{i}\right) =\left( b_{i}-1\right)
\delta _{i}.  \label{ab}
\end{equation}%
If $\delta _{i}\neq 0,1$ we determine one of $\left( a_{i},b_{i}\right) $ in
terms of the other through (\ref{ab}). We define $a_{i}$ in terms of $b_{i}$
for $\delta _{i}<.5$ and $b_{i}$ in terms of $a_{i}$ for $\delta _{i}>.5$;
this ensures that (\ref{phipoly}) below is symmetric. If $t_{1}=-1$ then $%
\delta _{1}=0$ and we set $a_{1}=1$. If $t_{p}=1$ then $\delta _{p}=1$ and
we set $b_{p}=1$. In each case the remaining parameter is set equal to $%
1/c $, so that the density tends to a point mass at $t_{i}$ as $\nu
\rightarrow 0$ and to uniformity as $\nu \rightarrow 1$.

The final density $\phi \left( \cdot \right) $ from which the design points
are to be sampled is a weighted average of those at (\ref{density}), with
weights proportional to the lengths $\left\vert I_{i}\right\vert $ of the $%
I_{i}$. Since $\left\vert J_{i}\right\vert =c\left\vert I_{i}\right\vert $
we obtain 
\begin{equation}
\phi \left( x;\nu \right) =\frac{1}{2c}\sum_{i=1}^{p}\beta
_{a_{i},b_{i}}\left( \frac{x-k_{i}}{\left\vert J_{i}\right\vert }\right)
I\left( x\in J_{i}\right) .  \label{phipoly}
\end{equation}%
Motivated by the designs of \S \ref{section: SLR} we recommend stratified
sampling, by which the sample consists of $\approx n\left\vert
I_{i}\right\vert /2$ points drawn from (\ref{density}), subject to an
appropriate rounding procedure.

\subsection{Polynomial regression\label{section: example2}}

We illustrate these proposals in the context of approximate polynomial
responses of degrees $p-1=1,2,3$. As also suggested in `Heuristic 5.1' of
Waite and Woods (2022, p.\ 1462), $t^{\ast }$ will consist of the support
points of the classical I-optimal designs. These I-optimal designs $\xi
^{\ast }$are obtained from Lemma 3.2 of Studden (1977), and are as
follows.\smallskip

$p=2$: $\xi ^{\ast }\left( \pm 1\right) =.5$,

$p=3$: $\xi ^{\ast }\left( \pm 1\right) =.25$, $\xi ^{\ast }\left( 0\right)
=.5$,

$p=4$: $\xi ^{\ast }\left( \pm 1\right) =\frac{1}{2\left( 1+\sqrt{5}\right) }%
\approx .1545$, $\xi ^{\ast }\left( \pm \frac{1}{\sqrt{5}}\approx \pm
.4472\right) =\frac{\sqrt{5}}{2\left( 1+\sqrt{5}\right) }\approx .3455$.
\smallskip

\begin{figure}[ptb]
\centering
\includegraphics[scale=1]{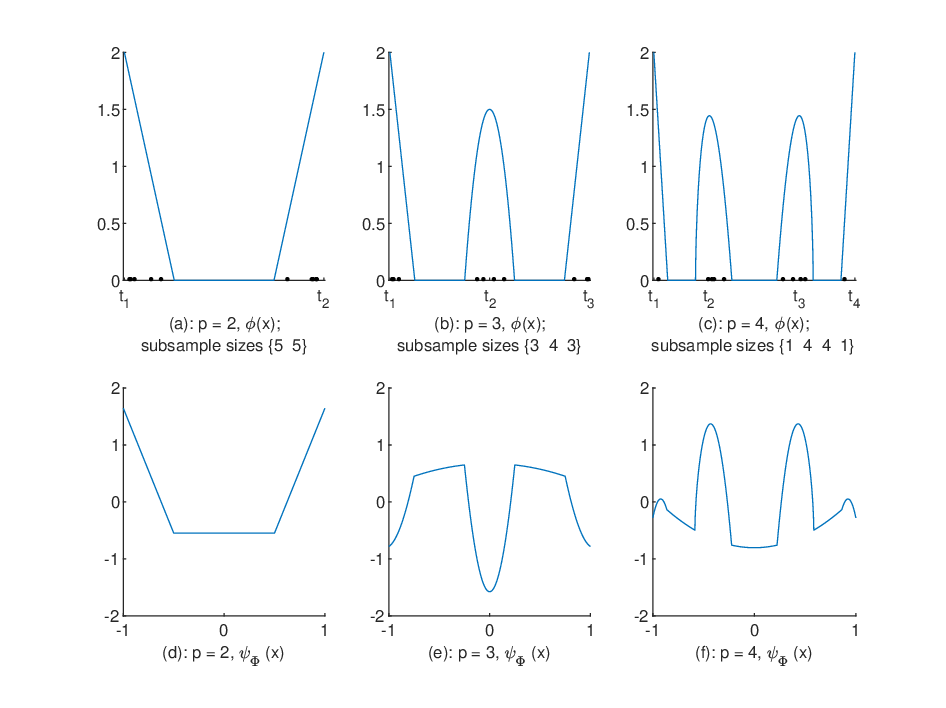}
\caption{Top: Cluster design densities $%
\protect\phi _{n}\left( x;\protect\nu =.5\right) $; typical stratified
samples using weights (a) $\left\{ .5,.5\right\} $, (b) $\left\{
.25,.5,.25\right\} $, (c) $\left\{ .14,.36,.36,.14\right\} $. Bottom: Scaled
least favourable contaminants $\protect\sqrt{n}\protect\psi _{\Phi }\left(
x\right) /\protect\tau $.}
\label{fig:cluster1}
\end{figure}

Figure \ref{fig:cluster1} gives the sampling densities (\ref{phipoly}),
together with the subsample sizes when $n=10$. Figure \ref{fig:cluster1}(a)
gives output for the approximate linear model, with a maximum \textsc{imse},
as at (\ref{max loss}),\ of $I_{\nu }\left( \Phi \right) =$ $2.804$. This
compares very favourably with the design of Figure \ref{fig:jitter1},
especially given that its construction does not require the minimax design
to be given. This latter point is especially germane for the design of
Figure \ref{fig:cluster1}(b), since it is the analogue of the absolutely
continuous minimax designs for approximate quadratic regression derived --
with substantial theoretical and computational difficulty -- by Shi et al.\ (2003) using methods of non-smooth optimization and by Daemi and Wiens
(2013) using completely different methods.

\begin{figure}[ptb]
\centering
\includegraphics[scale=1]{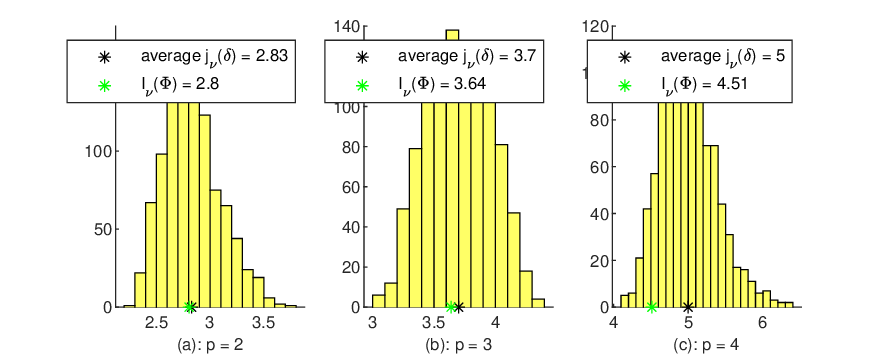}
\caption{Values of $j_{%
\protect\nu }\left( \protect\delta \right) $ from $1000$ cluster designs for
polynomial regression ($\protect\nu =.5$, $n=10$) and their averages,
estimating $J_{\protect\nu }\left( \Phi \right) $.}
\label{fig:cluster2}
\end{figure}

Figure \ref{fig:cluster2} gives values of $j_{\nu }\left( \delta
\right) $ from $1000$ random designs, together with their average,
estimating $J_{\nu }\left( \Phi \right) $. On average the random designs
perform almost as well against $\psi _{\Phi }$ -- plotted in Figure \ref%
{fig:cluster1}(d)-(f) -- as the continuous design $\Phi $.

It is interesting to note -- especially for the design of Figure \ref%
{fig:cluster1}(c) -- the close agreement between the I-optimal design
weights above, and the weights used in the computation of $\phi$  and detailed in the caption.

See Table 1, where the variance and maximum squared bias components of $%
I_{\nu }$ are presented for the designs of Figure \ref{fig:cluster1} ($\nu
=.5$) and for the corresponding designs with $\nu =.04$, very closely
approximating the I-optimal design ($\nu =0$) with maximum loss $%
I_{0}=\infty $. That the robustness of the cluster designs is achieved for
such a modest premium in terms of increased variance is both startling and
encouraging.

\begin{center}
\begin{table}[tb] \centering%
\begin{tabular}{ccccccccc}
\multicolumn{9}{c}{Table 1. \ Performance measures for the designs of Figure %
\ref{fig:cluster1}.} \\ \hline
& \multicolumn{2}{c}{variance} &  & \multicolumn{2}{c}{max sqd. bias} &  & 
\multicolumn{2}{c}{$I_{\nu }$} \\ \cline{2-3}\cline{5-6}\cline{8-9}
& $\nu =.5$ & $\nu =.04$ &  & $\nu =.5$ & $\nu =.04$ &  & $\nu =.5$ & $\nu
=.04$ \\ \cline{2-3}\cline{5-6}\cline{8-9}
$p=2$ & $2.94$ & $2.67$ &  & $2.67$ & $319$ &  & $2.80$ & $15.3$ \\ 
$p=3$ & $4.65$ & $4.27$ &  & $2.62$ & $213$ &  & $3.64$ & $12.6$ \\ 
$p=4$ & $6.49$ & $6.02$ &  & $2.54$ & $193$ &  & $4.51$ & $13.5$ \\ 
\hline\hline
\end{tabular}%
\end{table}%
\end{center}

\section{Multidimensional cluster designs\label{section: clusteringkd}}

\begin{figure}[tbp]
\centering
\includegraphics[scale=.75]{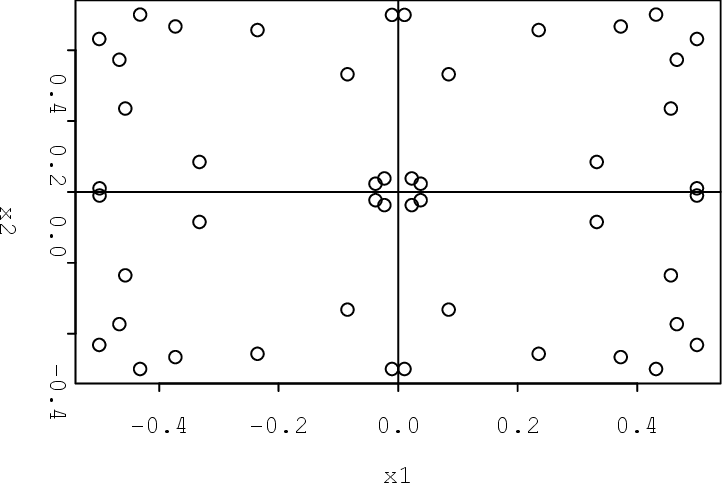}
\caption{Design for fitting a full second
order model; $n=48$.}
\label{fig:heo et al}
\end{figure}

See Figure \ref{fig:heo et al}, where a robust design, derived for fitting a
full second order bivariate model -- intercept, linear, quadratic and
interaction terms -- is depicted. It is a discrete implementation of a
design density, optimally robust against model misspecifications in a
certain parametric class of densities - see Heo et al.\ (2001) for details.
This design can roughly be described as an inscribed Central Composite
Design (CCD) with `clustering' in place of replication. It serves as
motivation for the ideas of this section, which we illustrate in the context
of $k$-dimensional, spherical CCD designs as are often used to fit second
order models. Such designs utilize $2^{k}+2k+1$ points $\left\{ \boldsymbol{t%
}_{i}\right\} $ consisting of $2^{k}$ corner points with $\boldsymbol{t}%
_{i,j}=\pm 1$ ($j=1,...,k$), $2k$ axial points $\boldsymbol{t}_{i}=\left(
0,...,\pm \sqrt{k},...,0\right) $ and a centre point $\boldsymbol{t}%
_{i}=\left( 0,...,0,...,0\right) $.

In this and other multidimensional cases we propose choosing design points
from spherical densities concentrated on neighbourhoods of the $\boldsymbol{t%
}_{i}$. A spherical density on a $k$-dimensional hypersphere%
\begin{equation*}
\mathcal{S}^{\left( k\right) }\left( \boldsymbol{t},R\right) =\left\{ 
\boldsymbol{x}\left\vert \left\Vert \boldsymbol{x}-\boldsymbol{t}\right\Vert
\leq R\right. \right\}
\end{equation*}%
with centre $\boldsymbol{t}$ and radius $R$, in which the scaled norm $%
\left\Vert \boldsymbol{x}-\boldsymbol{t}\right\Vert /R$ has a $Beta\left(
k,b\right) $ density, is given by 
\begin{equation*}
f^{(k)}\left( \boldsymbol{x};\boldsymbol{t},R,b\right) =\frac{\Gamma \left( 
\frac{k}{2}\right) }{2\pi ^{k/2}R^{k}\beta \left( k,b\right) }\cdot \left( 1-%
\frac{\left\Vert \boldsymbol{x}-\boldsymbol{t}\right\Vert }{R}\right)
^{b-1}I\left( \boldsymbol{x}\in \mathcal{S}^{\left( k\right) }\left( 
\boldsymbol{t},R\right) \right) .
\end{equation*}%
Such a density has mode $\boldsymbol{t}$ and approaches a point mass at $%
\boldsymbol{t}$ as $b\rightarrow \infty $, and uniformity as $b\rightarrow 1$%
. The choice of $k$ as the first parameter of the beta density ensures that $%
f$ is decreasing in $\left\Vert \boldsymbol{x}-\boldsymbol{t}\right\Vert $
and square integrable (required for the evaluation of the matrix $%
\boldsymbol{K}_{\Phi }$ as at (\ref{max loss})).

A sample value $\boldsymbol{x}$ from $f^{(k)}\left( \boldsymbol{x};%
\boldsymbol{t},R,b\right) $ is $\boldsymbol{x=t}+R\boldsymbol{y}$, where $%
\boldsymbol{y}\sim f^{(k)}\left( \boldsymbol{\cdot };\boldsymbol{0}%
,1,b\right) $ obtained by drawing a value of $\rho =\left\Vert \boldsymbol{y}%
\right\Vert \sim Beta(k,b)$ and, independently, drawing angles $\theta _{i}$%
, $-\pi /2<\theta _{i}\leq \pi /2$ ($i=1,...,k-2$) \ with densities $\psi
_{i}\left( \theta \right) =\cos ^{k-i-1}\theta \left/ \beta \left(
1/2,\left( k-i\right) /2\right) \right. $ -- equivalently, $\cos ^{2}\theta
_{i}$ $\sim Beta(\frac{1}{2},\frac{k-i}{2})$ -- and $\theta _{k-1}\sim
Unif\left( -\pi ,\pi \right) $. Then 
\begin{eqnarray*}
y_{1} &=&\rho \sin \theta _{1}, \\
y_{2} &=&\rho \cos \theta _{1}\sin \theta _{2}, \\
y_{3} &=&\rho \cos \theta _{1}\cos \theta _{2}\sin \theta _{3}, \\
&&\cdot \cdot \cdot \\
y_{k-1} &=&\rho \cos \theta _{1}\cos \theta _{2}\cdot \cdot \cdot \cos
\theta _{k-2}\sin \theta _{k-1}, \\
y_{k} &=&\rho \cos \theta _{1}\cos \theta _{2}\cdot \cdot \cdot \cos \theta
_{k-2}\cos \theta _{k-1}.
\end{eqnarray*}%
To sample $\theta _{i}$ for $i<k-1$ we draw $z\sim Beta(\frac{1}{2},\frac{k-i%
}{2})$ and set $\theta _{i}=\pm \arccos \sqrt{z}$, each with probability $%
1/2 $.

\subsection{Two dimensional cluster designs on tessellations \label{section:
tessellations}}

\begin{figure}[tb]
\centering
\includegraphics[trim={1.5cm 0 1cm 0},clip]{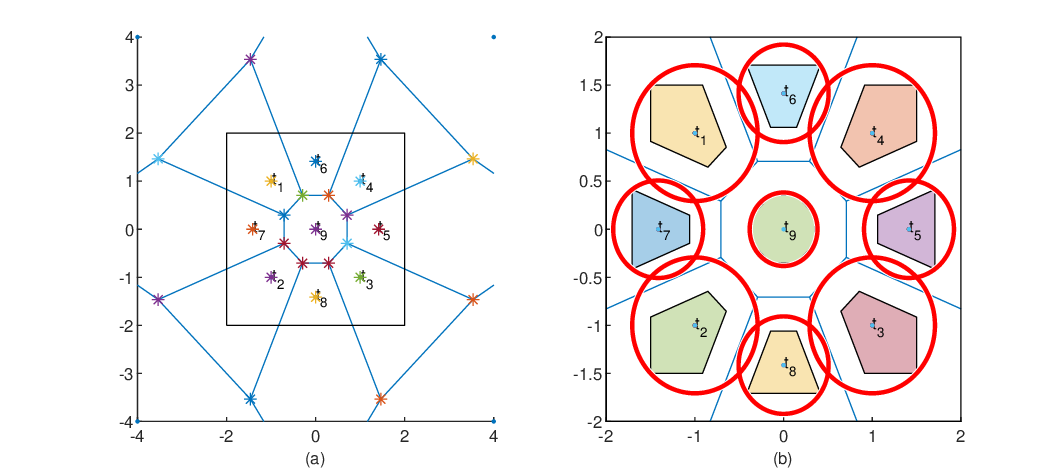}
\caption{(a) Voronoi tessellation
generated by the points $\left\{ \boldsymbol{t}_{i}\right\} $. (b)
Tessellation restricted to $\protect\chi =\left[ -2,2\right] ^{2}$ with
subtiles $\left\{ J_{i}\left( .25\right) \right\} $ and enclosing circles $%
\left\{ \mathcal{S}^{\left( 2\right) }\left( \boldsymbol{t}_{i},R_{i}\right)
\right\} $.}
\label{fig:tess1&2}
\end{figure}

\begin{figure}[tb]
\centering
\includegraphics[scale=1]{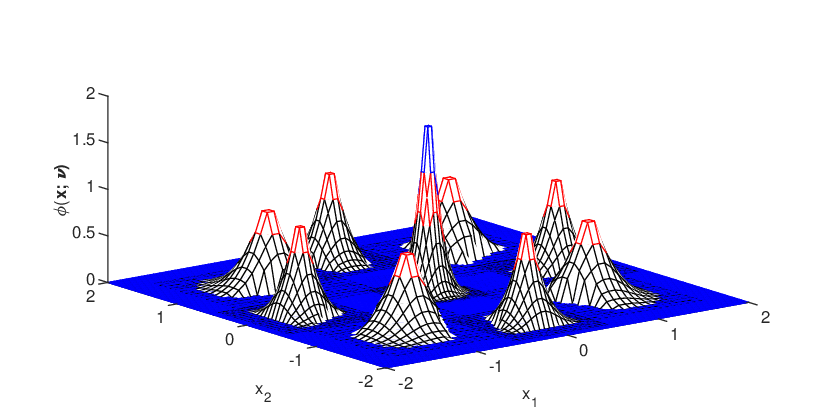}
\caption{Sampling density $\protect\phi %
\left( \boldsymbol{x};\protect\nu =.5\right) $ constructed for a robust,
clustered CCD in two dimensions.}
\label{fig:tess3}
\end{figure}

\begin{figure}[tbp]
\centering
\includegraphics[scale=1]{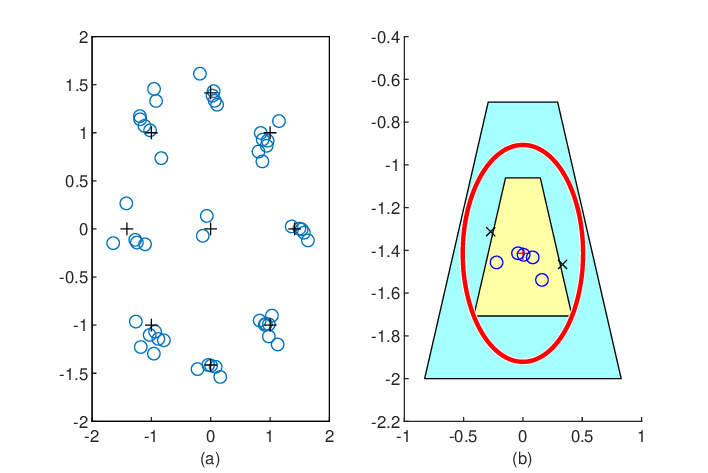}
\caption{(a) A typical random CCD of
size 50; $\protect\nu =.5$. (b) Details of the subsample of 5 points (`o')
drawn from $J_{8}\left( .25\right) $. Rejected points are marked as `x'.}
\label{fig:tess4}
\end{figure}

\begin{figure}[tbp]
\centering
\includegraphics[scale=1]{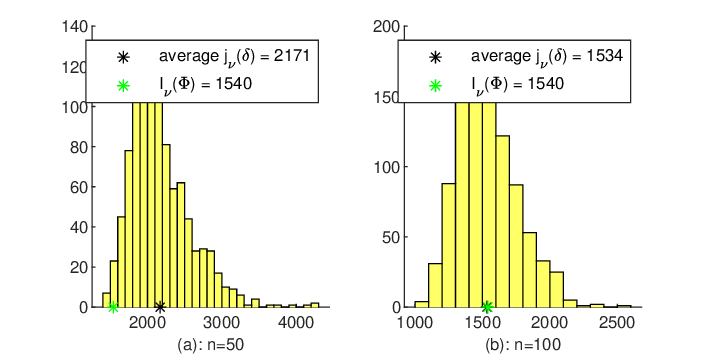}
\caption{Values of $j_{\protect\nu }\left( \protect%
\delta \right) $ from 1000 stratified, clustered CCDs ($\protect\nu =.5$)
and their averages, estimating $J_{\protect\nu }\left( \Phi \right) $. (a) $%
n=50$, (b) $n=100$.}
\label{fig:tess5}
\end{figure}

In Figure \ref{fig:tess1&2}(a), the nine points $\left\{ \boldsymbol{t}%
_{i}\right\} $ which are displayed consist of four corner points $(-1,\pm 1)$%
, $(1,\pm 1)$, four axial points $\left( \pm \sqrt{2},0\right) $, $\left(
0,\pm \sqrt{2}\right) $ and the centre point $\left( 0,0\right) $. These are
the generators of the \textit{Voronoi tessellation} pictured - a tiling with
the property that, within the tile $T_{i}$ containing $\boldsymbol{t}_{i}$,
all points are closer to $\boldsymbol{t}_{i}$ than to any $\boldsymbol{t}%
_{j},$ $j\neq i$. Figure \ref{fig:tess1&2}(b) gives a more detailed
depiction of the tessellation, restricted to the design space $\chi =\left[
-2,2\right] \times \lbrack -2,2]$. Within each tile $T_{i}$, of area $%
\left\vert T_{i}\right\vert $, we have also plotted a subtile $J_{i}\left(
c\right) $ which is a contraction of $T_{i}$ with fixed point $\boldsymbol{t}%
_{i}$ and area $\left\vert J_{i}\left( c\right) \right\vert =c\left\vert
T_{i}\right\vert $. These are then the analogues of the subintervals $%
J_{i}\left( c\right) \subseteq I_{i}$ from \S \ref{section: clustering1d},
and `$c$' has the same interpretation -- the fraction of the design space to
be sampled. Surrounding each $J_{i}\left( c\right) $ is the smallest
enclosing circle $\mathcal{S}^{\left( 2\right) }\left( \boldsymbol{t}%
_{i},R_{i}\left( c\right) \right) $.

We sample design points from $\mathcal{S}^{\left( 2\right) }\left( 
\boldsymbol{t}_{i},R_{i}\left( c\right) \right) $, accepting only those
points which lie in $J_{i}\left( c\right) $. We specify $b=1/c$ and $c=\nu
^{2}$, then $f^{(2)}\left( \boldsymbol{x};\boldsymbol{t}_{i},R_{i}\left(
c\right) ,b\right) $ approaches a point mass at $\boldsymbol{t}_{i}$ as $\nu
\rightarrow 0$, and uniformity on $\mathcal{S}^{\left( 2\right) }\left( 
\boldsymbol{t}_{i},R_{i}\left( c\right) \right) \supseteq T_{i}$ as $\nu
\rightarrow 1$. \ With 
\begin{equation*}
q_{i}(\nu )=\int_{J_{i}\left( c\right) }f^{(2)}\left( \boldsymbol{x};%
\boldsymbol{t}_{i},R_{i}\left( c\right) ,b\right) \mu \left( d\boldsymbol{x}%
\right) ,
\end{equation*}%
the density of those points accepted into the design upon being drawn from $%
\mathcal{S}^{\left( 2\right) }\left( \boldsymbol{t}_{i},R_{i}\left( c\right)
\right) $ is%
\begin{equation*}
\frac{f^{(2)}\left( \boldsymbol{x};\boldsymbol{t}_{i},R_{i}\left( c\right)
,b\right) }{q_{i}(\nu )}I\left( \boldsymbol{x}\in J_{i}\left( c\right)
\right) .
\end{equation*}%
We again do stratified sampling, with weights $\omega _{i}=\left\vert
T_{i}\right\vert \left/ \sum \left\vert T_{i}\right\vert \right. $
proportional to the area $\left\vert T_{i}\right\vert $, whence the density
of the design on $\chi $ is%
\begin{equation*}
\phi \left( \boldsymbol{x};\nu \right) =\sum_{i=1}^{9}\frac{\omega _{i}}{%
q_{i}(\nu )}f^{(2)}\left( \boldsymbol{x};\boldsymbol{t}_{i},R_{i}\left(
c\right) ,b\right) I\left( \boldsymbol{x}\in J_{i}\left( c\right) \right) .
\end{equation*}%
See Figure \ref{fig:tess3}. Although we evaluate $q_{i}(\nu )$ by numerical
integration, an estimate can be computed after the sampling is done; it is
the proportion of those points which were drawn from $\mathcal{S}^{\left(
2\right) }\left( \boldsymbol{t}_{i},R_{i}\left( c\right) \right) $ and then
accepted into the sample. This estimate turns out to be quite accurate if an
artificially large sample is simulated.

Figure \ref{fig:tess4} illustrates the results of applying the methods of
the preceding discussion. We chose a total sample size of $n=50$, $\nu =.5$,
and obtained subsample sizes $n_{i}=n\omega _{i}$, rounded to $\left\{
7,7,7,7,5,5,5,5,2\right\} $ with each corner point being allocated $7$, each
axial point being allocated $5$, and the remaining $2$ in the centre. The
entire sample is shown in Figure \ref{fig:tess4}(a), with Figure \ref%
{fig:tess4}(b) illustrating the details for Tile $8$. The required $5$
points were found after $7$ points were drawn from $\mathcal{S}_{8}$. In
all, $13$ points were rejected as not belonging to the appropriate subtile.

We repeated this with a total sample size of $n=100$. See Figure \ref%
{fig:tess5}. With the larger sample size the random designs seem to more
accurately duplicate the behaviour of the parent design $\Phi $ -- a
phenomenon noticed as well in the one dimensional cluster designs of the
previous section.

\subsection{Extensions to $k>2$ \label{section: k>2}}

\begin{figure}[tb]
\centering
\includegraphics[trim={0 4cm 0 0},clip]{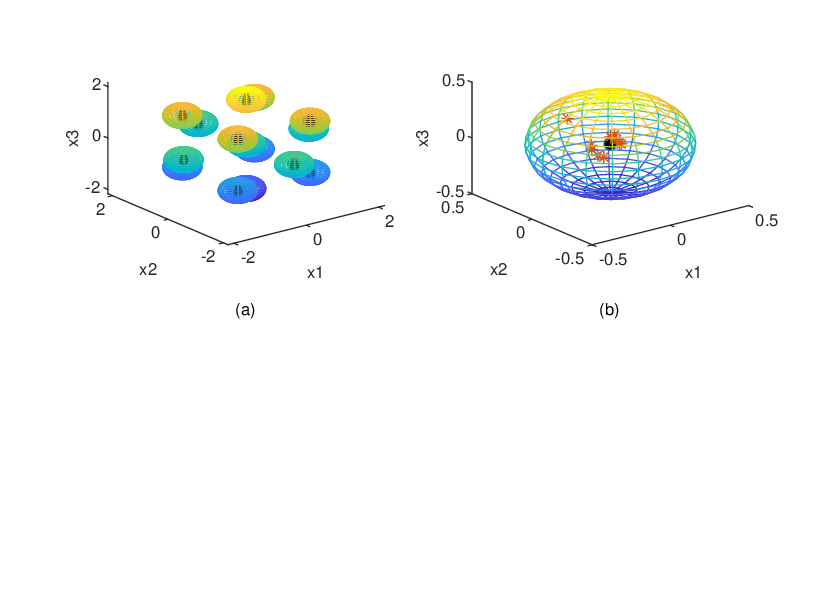}
\caption{(a) Spheres $\mathcal{S}%
^{\left( 3\right) }\left( \boldsymbol{t}_{i},r_{0}\protect\nu \right) $, $%
i=1,...,15$ for a three dimensional spherical CCD and $\protect\nu =.5$. (b)
Centre sphere $\mathcal{S}^{\left( 3\right) }\left( \boldsymbol{0},r_{0}%
\protect\nu \right) $ with sampled points.}
\label{fig:ccd3}
\end{figure}

Although the theory of \S \ref{section: tessellations} extends easily to
higher dimensions, the lack of appropriate software for constructing and
manipulating Voronoi tessellations becomes a severe drawback. But the
general idea of sampling from spherical distributions centred on small
neighbourhoods of the $\left\{ \boldsymbol{t}_{i}\right\} $ can still be
applied, albeit in a less structured manner. Let $\left\{ \boldsymbol{t}%
_{i}\right\} _{i=1}^{q}$ be the $q=2^{k}+2k+1$ support points of a spherical
CCD in variables $\boldsymbol{x}=\left( x_{1},...,x_{k}\right) ^{\prime }$%
, as described at the beginning of this section. The minimum distance
between these points is $\min \left( 2,\sqrt{k}\right) $, and so
hyperspheres $\mathcal{S}\left( \boldsymbol{t}_{i},r_{0}\right) $ centred at
the $\boldsymbol{t}_{i}$ and with radius $r_{0}=\min \left( 1,\sqrt{k}%
/2\right) $ are disjoint. \ Define subspheres 
\begin{equation*}
J_{i}\left( c\right) =\mathcal{S}^{\left( k\right) }\left( \boldsymbol{t}%
_{i},r_{0}c^{1/k}\right) ,0<c\leq 1.
\end{equation*}%
Then $\int I\left( \boldsymbol{x}\in J_{i}\left( c\right) \right) d%
\boldsymbol{x}=\left\vert J_{i}\left( c\right) \right\vert =c\left\vert 
\mathcal{S}\left( \boldsymbol{t}_{i},r_{0}\right) \right\vert $. The density
of $\boldsymbol{x}$ on $J_{i}\left( c\right) $ is $f^{(k)}\left( \boldsymbol{%
\cdot };\boldsymbol{t}_{i},r_{0}c^{1/k},b\right) $. We again specify $b=1/c$
and $c=\nu ^{k}$. Then for user chosen weights $\left\{ \omega _{i}\right\} $
the sampling density is%
\begin{equation*}
\phi \left( \boldsymbol{x};\nu \right) =\sum_{i=1}^{q}\omega
_{i}f^{(k)}\left( \boldsymbol{x};\boldsymbol{t}_{i},r_{0}\nu ,1/\nu \right)
I\left( \boldsymbol{x}\in \mathcal{S}^{\left( k\right) }\left( \boldsymbol{t}%
_{i},r_{0}\nu \right) \right) .
\end{equation*}

See Figure \ref{fig:ccd3} for an example with $k=3$. We sampled a design of
size $n=80$ with subsamples sizes $n_{i}=5$ ($i<15$) and $n_{15}=10$.

\section*{Appendix}

\setcounter{equation}{0} \setcounter{subsection}{0} \renewcommand{%
\theequation}{A.\arabic{equation}} \renewcommand{\thesubsection}{A.%
\arabic{subsection}}

\subsection{Derivations for \S \protect\ref{section: random designs}}

For an \ $n$-point design $D=\left\{ \boldsymbol{x}_{i}\right\} _{i=1}^{n}$
with design measure $\delta =n^{-1}\sum \delta _{\boldsymbol{x}_{i}}$ define%
\begin{equation*}
\boldsymbol{F}=\left( \boldsymbol{f}\left( \boldsymbol{x}_{1}\right) ,\cdot
\cdot \cdot ,\boldsymbol{f}\left( \boldsymbol{x}_{n}\right) \right) ^{\prime
},\boldsymbol{\psi }_{\Phi }=\left( \psi _{\Phi }\left( \boldsymbol{x}%
_{1}\right) ,...,\psi _{\Phi }\left( \boldsymbol{x}_{n}\right) \right)
^{\prime },\text{ }\boldsymbol{D}_{\phi }=diag\left( \phi \left( \boldsymbol{%
x}_{1}\right) ,...,\phi \left( \boldsymbol{x}_{n}\right) \right) .
\end{equation*}%
Then $\boldsymbol{M}_{\phi }=\frac{1}{n}\boldsymbol{F}^{\prime }\boldsymbol{D%
}_{\phi }\boldsymbol{F}$. Define as well 
\begin{eqnarray*}
\boldsymbol{M}_{\delta } &=&\int_{\chi }\boldsymbol{f}\left( \boldsymbol{x}%
\right) \boldsymbol{f}^{\prime }\left( \boldsymbol{x}\right) \delta \left( d%
\boldsymbol{x}\right) =\frac{1}{n}\sum_{\boldsymbol{x}_{i}\in D}\boldsymbol{f%
}\left( \boldsymbol{x}_{i}\right) \boldsymbol{f}^{\prime }\left( \boldsymbol{%
x}_{i}\right) =\frac{1}{n}\boldsymbol{F}^{\prime }\boldsymbol{F},\text{ \ }
\\
\boldsymbol{b}_{\psi _{\Phi },\delta } &=&\int_{\chi }\boldsymbol{f}\left( 
\boldsymbol{x}\right) \psi _{\Phi }\left( \boldsymbol{x}\right) \delta
\left( d\boldsymbol{x}\right) =\frac{1}{n}\sum_{\boldsymbol{x}_{i}\in D}%
\boldsymbol{f}\left( \boldsymbol{x}_{i}\right) \psi _{\Phi }\left( 
\boldsymbol{x}_{i}\right) =\frac{1}{n}\boldsymbol{F}^{\prime }\boldsymbol{%
\psi }_{\Phi }.
\end{eqnarray*}%
Using (\ref{least favourable}),

\begin{equation}
\psi _{\Phi }\left( \boldsymbol{x}_{i}\right) =\mathbf{r}_{\Phi }^{\prime
}\left( \boldsymbol{x}_{i}\right) \mathbf{\beta }_{\Phi }=\frac{\tau }{\sqrt{%
n}}\left( \phi \left( \boldsymbol{x}_{i}\right) \boldsymbol{f}^{\prime
}\left( \boldsymbol{x}_{i}\right) -\boldsymbol{f}^{\prime }\left( 
\boldsymbol{x}_{i}\right) \boldsymbol{A}^{-1}\boldsymbol{M}_{\Phi }\right) 
\boldsymbol{G}_{\Phi }^{-1/2}\mathbf{\beta }_{\Phi }\mathbf{,}
\label{psi(xi)}
\end{equation}%
so that%
\begin{eqnarray}
\boldsymbol{\psi }_{\Phi } &=&\frac{\tau }{\sqrt{n}}\left( \boldsymbol{D}%
_{\phi }\boldsymbol{F}\mathbf{-}\boldsymbol{FA}^{-1}\boldsymbol{M}_{\Phi
}\right) \boldsymbol{G}_{\Phi }^{-1/2}\mathbf{\beta }_{\Phi }\mathbf{,}
\label{psi_phi} \\
\boldsymbol{b}_{\psi _{\Phi },\delta } &=&\frac{\tau }{\sqrt{n}}\left( 
\boldsymbol{M}_{\phi }-\boldsymbol{M}_{\delta }\boldsymbol{A}^{-1}%
\boldsymbol{M}_{\Phi }\right) \boldsymbol{G}_{\Phi }^{-1/2}\mathbf{\beta }%
_{\Phi }\mathbf{.}  \label{b}
\end{eqnarray}%
From (\ref{imse}),%
\begin{equation*}
\text{\textsc{imse}}\left( \delta |\psi _{\Phi }\right) =\frac{\sigma
_{\varepsilon }^{2}}{n}tr\left( \boldsymbol{AM}_{\delta }^{-1}\right) +%
\boldsymbol{b}_{\psi _{\Phi },\delta }^{\boldsymbol{\prime }}\boldsymbol{M}%
_{\delta }^{-1}\boldsymbol{AM}_{\delta }^{-1}\boldsymbol{b}_{\psi _{\Phi
},\delta }+\int_{\chi }\psi _{\Phi }^{2}\left( x\right) \mu \left( dx\right)
;
\end{equation*}%
substituting (\ref{psi_phi}) and (\ref{b}) gives 
\begin{eqnarray*}
\text{\textsc{imse}}\left( \delta |\psi _{\Phi }\right) &=&\frac{\sigma
_{\varepsilon }^{2}}{n}tr\left( \boldsymbol{AM}_{\delta }^{-1}\right) +\frac{%
\tau ^{2}}{n}\gamma _{\delta },\text{ for} \\
\gamma _{\delta } &=&\mathbf{\beta }_{\Phi }^{\prime }\boldsymbol{G}_{\Phi
}^{-1/2}\left( \boldsymbol{M}_{\phi }\boldsymbol{M}_{\delta }^{-1}-%
\boldsymbol{M}_{\Phi }\boldsymbol{A}^{-1}\right) \boldsymbol{A}\left( 
\boldsymbol{M}_{\delta }^{-1}\boldsymbol{M}_{\phi }-\boldsymbol{A}^{-1}~%
\boldsymbol{M}_{\Phi }\right) \boldsymbol{G}_{\Phi }^{-1/2}\mathbf{\beta }%
_{\Phi }+1.
\end{eqnarray*}%
Now (\ref{E(imse)}) and (\ref{J}) are immediate.

\subsection{Derivations of (\protect\ref{easy J}) and (\protect\ref{PsiPhi})}

To evaluate (\ref{J}) and establish (\ref{easy J}) we first note that since $%
\phi \left( x_{i}\right) \equiv \left( 2c\right) ^{-1}$ on its support, we
have that $\boldsymbol{M}_{\phi }\boldsymbol{M}_{\delta }^{-1}=\left(
2c\right) ^{-1}\boldsymbol{I}_{2}$, and then (since $\mathbf{\beta }_{\Phi
}^{\prime }\mathbf{\beta }_{\Phi }=1$) 
\begin{eqnarray}
\gamma _{\delta } &=&\mathbf{\beta }_{\Phi }^{\prime }\boldsymbol{G}_{\Phi
}^{-1/2}\left( \frac{1}{2c}\boldsymbol{I}_{2}-\boldsymbol{M}_{\Phi }%
\boldsymbol{A}^{-1}\right) \boldsymbol{A}\left( \frac{1}{2c}\boldsymbol{I}%
_{2}-\boldsymbol{A}^{-1}\boldsymbol{M}_{\Phi }\right) \boldsymbol{G}_{\Phi
}^{-1/2}\mathbf{\beta }_{\Phi }+1  \notag \\
&=&\mathbf{\beta }_{\Phi }^{\prime }\boldsymbol{G}_{\Phi }^{-1/2}\left( 
\frac{1}{4c^{2}}\boldsymbol{A}-\frac{1}{c}\boldsymbol{M}_{\Phi }+\boldsymbol{%
M}_{\Phi }\boldsymbol{A}^{-1}\boldsymbol{M}_{\Phi }\right) \boldsymbol{G}%
_{\Phi }^{-1/2}\mathbf{\beta }_{\Phi }+1  \notag \\
&=&\mathbf{\beta }_{\Phi }^{\prime }\boldsymbol{G}_{\Phi }^{-1/2}\left( 
\frac{1}{4c^{2}}\boldsymbol{A}-\frac{1}{c}\boldsymbol{M}_{\Phi }+\boldsymbol{%
H}_{\Phi }+\boldsymbol{G}_{\Phi }\right) \boldsymbol{G}_{\Phi }^{-1/2}%
\mathbf{\beta }_{\Phi }.  \label{gamma_new}
\end{eqnarray}%
\bigskip A calculation gives 
\begin{equation}
\boldsymbol{G}_{\Phi }^{1/2}\boldsymbol{H}_{\Phi }^{-1}\boldsymbol{G}_{\Phi
}^{1/2}+\boldsymbol{I}_{2}=\frac{1}{c}diag\left( 1,\frac{1}{3\lambda
_{2}\left( c\right) }\right) ,  \label{diag}
\end{equation}%
so that the maximum eigenvalue is $\gamma _{0}$ and 
\begin{equation*}
\mathbf{\beta }_{\Phi }=\left\{ 
\begin{array}{cc}
\left( 1,0\right) ^{\prime }, & \text{if }\lambda _{2}\left( c\right) \geq
1/3, \\ 
\left( 0,1\right) ^{\prime }, & \text{if }\lambda _{2}\left( c\right) <1/3.%
\end{array}%
\right.
\end{equation*}%
The choice of $\mathbf{\beta }_{\Phi }$ is somewhat arbitrary if $\lambda
_{2}\left( c\right) =1/3$, since then (\ref{diag}) is a multiple of $%
\boldsymbol{I}_{2}$. \ We claim that%
\begin{equation}
\gamma _{\delta }\equiv \gamma _{0},  \label{claim}
\end{equation}%
from which (\ref{easy J}) follows, since then $\gamma _{\delta }$ does not
depend on the design and so is non-random.

To establish (\ref{claim}), use $\boldsymbol{A}=\boldsymbol{M}_{\Phi }%
\boldsymbol{H}_{\Phi }^{-1}\boldsymbol{M}_{\Phi }=4c^{2}\boldsymbol{K}_{\Phi
}\boldsymbol{H}_{\Phi }^{-1}\boldsymbol{K}_{\Phi }$ and $\boldsymbol{M}%
_{\Phi }=2c\boldsymbol{K}_{\Phi }$ in (\ref{gamma_new}) to obtain 
\begin{equation*}
\gamma _{\delta }=\mathbf{\beta }_{\Phi }^{\prime }\left[ \boldsymbol{G}%
_{\Phi }^{-1/2}\left( \boldsymbol{K}_{\Phi }\boldsymbol{H}_{\Phi }^{-1}%
\boldsymbol{K}_{\Phi }-\boldsymbol{K}_{\Phi }\right) \boldsymbol{G}_{\Phi
}^{-1/2}\right] \mathbf{\beta }_{\Phi }.
\end{equation*}%
Substituting $\boldsymbol{K}_{\Phi }=\boldsymbol{G}_{\Phi }+\boldsymbol{H}%
_{\Phi }$, this becomes%
\begin{equation*}
\gamma _{\delta }=\mathbf{\beta }_{\Phi }^{\prime }\left[ \boldsymbol{G}%
_{\Phi }^{1/2}\boldsymbol{H}_{\Phi }^{-1}\boldsymbol{G}_{\Phi }^{1/2}+%
\boldsymbol{I}_{2}\right] \mathbf{\beta }_{\Phi }=\gamma _{0},
\end{equation*}%
as required.

\bigskip An evaluation of (\ref{psi(xi)}), using 
\begin{eqnarray*}
\boldsymbol{A} &=&diag\left( 2,2/3\right) ,\text{ }\boldsymbol{M}_{\Phi
}=diag\left( 1,\lambda _{2}\left( c\right) \right) ,\text{ and } \\
\boldsymbol{G}_{\Phi } &=&\frac{1}{2c}diag\left( \left( 1-c\right) ,\lambda
_{2}\left( c\right) \left( 1-3c\lambda _{2}\left( c\right) \right) \right) ,
\end{eqnarray*}%
gives 
\begin{equation*}
\psi _{\Phi }\left( x;c\right) =\frac{\tau }{\sqrt{n}}\left( \frac{2c\phi
_{n}\left( x;c\right) -c}{\sqrt{2c\left( 1-c\right) }},x\frac{2c\phi
_{n}\left( x;c\right) -3c\lambda _{2}\left( c\right) }{\sqrt{2c\lambda
_{2}\left( c\right) \left( 1-3c\lambda _{2}\left( c\right) \right) }}\right) 
\mathbf{\beta }_{\Phi }.
\end{equation*}%
Using $2c\phi _{n}\left( x;c\right) =I_{S}\left( x;c\right) $, and 
\begin{equation*}
\frac{c}{\gamma _{0}}=\left\{ 
\begin{array}{cc}
c, & \text{if }\lambda _{2}\left( c\right) \geq 1/3, \\ 
3c\lambda _{2}\left( c\right) & \text{if }\lambda _{2}\left( c\right) <1/3,%
\end{array}%
\right.
\end{equation*}%
this becomes%
\begin{equation*}
\psi _{\Phi }\left( x;c\right) =\frac{\tau }{\sqrt{n}}\cdot \left\{ 
\begin{array}{cc}
\frac{I_{S}\left( x;c\right) -\frac{1}{\gamma _{0}}}{\sqrt{2c\left( 1-\frac{1%
}{\gamma _{0}}\right) }}, & \text{if }\lambda _{2}\left( c\right) \geq 1/3,
\\ 
\frac{x}{\sqrt{\lambda _{2}}}\frac{I_{S}\left( x;c\right) -\frac{1}{\gamma
_{0}}}{\sqrt{2c\left( 1-\frac{1}{\gamma _{0}}\right) }} & \text{if }\lambda
_{2}\left( c\right) <1/3,%
\end{array}%
\right. ,
\end{equation*}%
which is (\ref{PsiPhi}).

\section{Acknowledgements}

This work was carried out with the support of the Natural Sciences and
Engineering Research Council of Canada. It has benefited from conversations
with Timothy Waite, University of Manchester and Xiaojian Xu, Brock
University. We are grateful for the incisive and helpful comments of the
reviewers.

\section*{References}

\begin{description}

\item Atkinson, A. C.\ (1996), \textquotedblleft The Usefulness of Optimum
Experimental Designs,\textquotedblright\ \textit{Journal of the Royal
Statistical Society (Series B)}; 58, 59-76.

\item Bischoff, W.\ (2010), \textquotedblleft An Improvement in the
Lack-of-Fit Optimality of the (Absolutely) Continuous Uniform Design in
Respect of Exact Designs,\textquotedblright\ in\ \textit{mODa 9 - Advances
in Model-Oriented Design and Analysis}, eds.\ Giovagnoli, G., Atkinson, A.
and Torsney, B.

\item Box, G.E.P., and Draper, N.R.\ (1959), \textquotedblleft A Basis for
the Selection of a Response Surface Design,\textquotedblright\ \textit{%
Journal of the American Statistical Association}, 54, 622-654.

\item Daemi, M., and Wiens, D.P.\ (2013), \textquotedblleft Techniques for
the Construction of Robust Regression Designs,\textquotedblright\ \textit{%
The Canadian Journal of Statistics}, 41, 679-695.

\item Fang, K.T.\ \& Wang, Y.\ (1994). \textit{Number-Theoretic Methods in
Statistics}. Chapman and Hall, London and New York.

\item Fang, Z., and Wiens, D.P.\ (2000), \textquotedblleft Integer-Valued,
Minimax Robust Designs for Estimation and Extrapolation in Heteroscedastic,
Approximately Linear Models,\textquotedblright\ \textit{Journal of the
American Statistical Association}, 95, 807-818.

\item Heo, G., Schmuland, B., and Wiens, D.P.\ (2001), \textquotedblleft
Restricted Minimax Robust Designs for Misspecified Regression
Models,\textquotedblright\ \textit{The Canadian Journal of Statistics}, 29,
117-128.

\item Huber, P.J.\ (1975), \textquotedblleft Robustness and
Designs,\textquotedblright\ in: \textit{A Survey of Statistical Design and
Linear Models}, ed.\ J.N.\ Srivastava, Amsterdam: North Holland, pp.\
287-303.

\item Shi, P., Ye, J., and Zhou, J. (2003), \textquotedblleft Minimax Robust
Designs for Misspecified Regression Models,\textquotedblright\ \textit{The
Canadian Journal of Statistics}, 31, 397-414.

\item Studden, W.J.\ (1977), \textquotedblleft Optimal Designs for
Integrated Variance in Polynomial Regression,\textquotedblright\ \textit{%
Statistical Decision Theory and Related Topics II}, ed.\ Gupta, S.S. and
Moore, D.S.\thinspace New York: Academic Press, pp.\ 411-420.

\item Waite, T.W., and Woods, D.C.\ (2022),\textquotedblleft Minimax
Efficient Random Experimental Design Strategies With Application to
Model-Robust Design for Prediction,\textquotedblright\ \textit{Journal of
the American Statistical Association}, 117, 1452-1465.

\item Wiens, D.P.\ (1990), \textquotedblleft Robust, Minimax Designs for
Multiple Linear Regression,\textquotedblright\ \textit{Linear Algebra and
Its Applications, Second Special Issue on Linear Algebra and Statistics};
127, 327 - 340.

\item Wiens, D.P.\ (1992), \textquotedblleft Minimax Designs for
Approximately Linear Regression,\textquotedblright\ \textit{Journal of
Statistical Planning and Inference,} 31, 353-371.

\item Wiens, D.P.\ (2015), \textquotedblleft Robustness of
Design\textquotedblright , Chapter 20, \textit{Handbook of Design and
Analysis of Experiments}, Chapman \& Hall/CRC.

\item Wiens, D.P.\ (2018), \textquotedblleft I-Robust and D-Robust Designs
on a Finite Design Space,\textquotedblright\ \textit{Statistics and Computing%
}, 28, 241-258.

\item Xu, X., and Yuen, W.K.\ (2011), \textquotedblleft Applications and
Implementations of Continuous Robust Designs,\textquotedblright\ \textit{%
Communications in Statistics - Theory and Methods}, 40, 969-988.
\end{description}

\end{document}